\title{Privileged users in zero-error transmission over a noisy channel}
\author{{Noga Alon\thanks{ Schools of Mathematics and Computer Science,
Raymond and Beverly Sackler Faculty of Exact Sciences, Tel Aviv
University, Tel Aviv, 69978, Israel. Email: nogaa@tau.ac.il.
Research supported in part by a USA-Israeli BSF grant, by the Israel
Science Foundation and by the Hermann Minkowski Minerva Center for
Geometry at Tel Aviv University.}} \quad {Eyal Lubetzky
\thanks{ School of Computer Science, Raymond and Beverly
Sackler Faculty of Exact Sciences, Tel Aviv University, Tel Aviv,
69978, Israel. Email: lubetzky@tau.ac.il. Research partially
supported by a Charles Clore Foundation Fellowship.}}}

\documentclass [11pt]{article}
\usepackage[]{graphicx,epsf,epsfig,amsmath,amssymb,amsfonts,latexsym,amsthm}

\newtheorem{theorem}{Theorem}[section]
\newtheorem{lemma}[theorem]{Lemma}

\newtheorem*{definition}{Definition}

\newtheorem{proposition}[theorem]{Proposition}

\renewcommand{\epsilon}{\varepsilon}

\newtheoremstyle{upright}%
        {8pt plus2pt minus4pt}%
        {8pt plus2pt minus4pt}%
        {\upshape}%
        {}%
        {\bfseries\scshape}%
        {:}%
        {1em}%
        {}%
\theoremstyle{upright}

\newcommand{\ignore}[1]{}

\begin{document}
\maketitle

\begin{abstract} The $k$-th power of a graph $G$ is the graph whose
vertex set is $V(G)^k$, where two distinct $k$-tuples are adjacent
iff they are equal or adjacent in $G$ in each coordinate. The
Shannon capacity of $G$, $c(G)$, is
$\lim_{k\to\infty}\alpha(G^k)^{\frac{1}{k}}$, where $\alpha(G)$
denotes the independence number of $G$. When $G$ is the
characteristic graph of a channel $\mathcal{C}$, $c(G)$ measures the
effective alphabet size of $\mathcal{C}$ in a zero-error protocol. A
sum of channels, $\mathcal{C}=\sum_i \mathcal{C}_i$, describes a
setting when there are $t\geq 2$ senders, each with his own channel
$\mathcal{C}_i$, and each letter in a word can be selected from
either of the channels. This corresponds to a disjoint union of the
characteristic graphs, $G=\sum_i G_i$. It is well known that $c(G)
\geq \sum_i c(G_i)$, and in \cite{NogaUnion} it is shown that in
fact $c(G)$ can be larger than any fixed power of the above sum.

We extend the ideas of \cite{NogaUnion} and show that for every
$\mathcal{F}$, a family of subsets of $[t]$, it is possible to
assign a channel $\mathcal{C}_i$ to each sender $i\in[t]$, such that
the capacity of a group of senders $X \subset [t]$ is high iff $X$
contains some $F \in \mathcal{F}$. This corresponds to a case where
only privileged subsets of senders are allowed to transmit in a high
rate. For instance, as an analogue to secret sharing, it is possible
to ensure that whenever at least $k$ senders combine their channels,
they obtain a high capacity, however every group of $k-1$ senders
has a low capacity (and yet is not totally denied of service). In
the process, we obtain an explicit Ramsey construction of an
edge-coloring of the complete graph on $n$ vertices by $t$ colors,
where every induced subgraph on $\exp\left(\Omega(\sqrt{\log n
\log\log n})\right)$ vertices contains all $t$ colors.
\end{abstract}

\section{Introduction}\label{sec::intro}
A channel $\mathcal{C}$ on an input alphabet $V$ and an output
alphabet $U$ maps each $x \in V$ to some $S(x)\subset U$, such that
transmitting $x$ results in one of the letters of $S(x)$. The
characteristic graph of the channel $\mathcal{C}$,
$G=G(\mathcal{C})$, has a vertex set $V$, and two vertices $x\neq
y\in V$ are adjacent iff $S(x)\cap S(y)\neq \emptyset$, i.e., the
corresponding input letters are confusable in the channel. Clearly,
a maximum set of predefined letters which can be transmitted in
$\mathcal{C}$ without possibility of error corresponds to a maximum
independent set in the graph $G$, and has cardinality $\alpha(G)$
(the independence number of $G$).

The strong product of two graphs, $G_1=(V_1,E_1)$ and
$G_2=(V_2,E_2)$ is the graph, $G_1 \cdot G_2$, on the vertex set
$V_1\times V_2$, where two distinct vertices
$(u_1,u_2)\neq(v_1,v_2)$ are adjacent iff for all $i=1,2$, either
$u_i=v_i$ or $u_i v_i \in E_i$. In other words, the pairs of
vertices in both coordinates are either equal or adjacent. This
product is associative and commutative, hence we can define $G^k$ to
be the $k$-th power of $G$, where two vertices $(u_1,\ldots,u_k)\neq
(v_1,\ldots,v_k)$ are adjacent iff for all $i=1,\ldots,k$, either
$u_i=v_i$ or $u_i v_i \in E(G)$.

Note that if $I,J$ are independent sets of two graphs, $G,H$, then
$I\times J$ is an independent set of $G \cdot H$. Therefore,
$\alpha(G^{n+m})\geq\alpha(G^n)\alpha(G^m)$ for every $m,n\geq 1$,
and by Fekete's lemma (cf., e.g., \cite{VanLintWilson}, p. 85), the
limit $\lim_{n\to\infty}\alpha(G^n)^{\frac{1}{n}}$ exists and equals
$\sup_n \alpha(G^n)^{\frac{1}{n}}$. This parameter, introduced by
Shannon in \cite{Shannon}, is the Shannon capacity of $G$, denoted
by $c(G)$.

When sending $k$-letter words in the channel $\mathcal{C}$, two
words are confusable iff the pairs of letters in each of their
$k$-coordinates are confusable. Thus, the maximal number of
$k$-letter words which can be sent in $\mathcal{C}$ without
possibility of error is precisely $\alpha(G^k)$, where
$G=G(\mathcal{C})$. It follows that for sufficiently large values of
$k$, the maximal number of $k$-letter words which can be sent
without possibility of error is roughly $c(G)^k$. Hence, $c(G)$
represents the effective alphabet size of the channel in zero-error
transmission.

The sum of two channels, $\mathcal{C}_1+\mathcal{C}_2$, describes
the setting where each letter can be sent from either of the two
channels, and letters from $\mathcal{C}_1$ cannot be confused with
letters from $\mathcal{C}_2$. The characteristic graph in this case
is the disjoint union $G_1+G_2$, where $G_i$ is the characteristic
graph of $\mathcal{C}_i$. Shannon showed in \cite{Shannon} that
$c(G_1+G_2)\geq c(G_1)+c(G_2)$ for every two graphs $G_1$ and $G_2$,
and conjectured that in fact $c(G_1+G_2)=c(G_1)+c(G_2)$ for all
$G_1$ and $G_2$. This was disproved in \cite{NogaUnion}, where the
author gives an explicit construction of two graphs $G_1,G_2$ with a
capacity $c(G_i)\leq k$, satisfying $c(G_1+G_2)\geq
k^{\Omega(\frac{\log k}{\log\log k})}$.

We extend the ideas of \cite{NogaUnion} and show that it is possible
to construct $t$ graphs, $G_i$ ($i\in[t]=\{1,2,\ldots ,t\}$), such
that for every subset $X\subseteq [t]$, the Shannon capacity of
$\sum_{i\in X}G_i$ is high iff $X$ contains some subset of a
predefined family $\mathcal{F}$ of subsets of $[t]$. This
corresponds to assigning $t$ channels to $t$ senders, such that
designated groups of senders $F\in\mathcal{F}$ can obtain a high
capacity by combining their channels ($\sum_{i\in F}\mathcal{C}_i$),
and yet every group of senders $X\subseteq[t]$ not containing any
$F\in\mathcal{F}$ has a low capacity. In particular, a choice of
$\mathcal{F}=\{F \subset [t]~:~|F|=k\}$ implies that every set $X$
of senders has a high Shannon capacity of $\sum_{i\in
X}\mathcal{C}_i$ if $|X|\geq k$, and a low capacity otherwise. The
following theorem, proved in Section \ref{sec-2}, formalizes the
claims above:
\begin{theorem}\label{thm-1} Let $T=\{1,\ldots,t\}$ for some fixed $t \geq 2$,
and let $\mathcal{F}$ be a family of subsets of $T$. For every
(large) $n$ it is possible to construct graphs $G_i$, $i\in T$, each
on $n$ vertices, such that the following two statements hold for all
$X \subseteq T$:
\begin{enumerate}
\item If $X$ contains some $F\in\mathcal{F}$, then $c(\sum_{i\in
X}G_i) \geq n^{1/|F|} \geq n^{1/t}$.
\item If $X$ does not contain any $F\in\mathcal{F}$, then $$c(\sum_{i\in
X}G_i) \leq \mathrm{e}^{(1+o(1))\sqrt{2\log n \log \log n}}~,$$
where the $o(1)$-term tends to 0 as $n\to\infty$.
\end{enumerate}
\end{theorem}
As a by-product, we obtain the following Ramsey construction, where
instead of forbidding monochromatic subgraphs, we require ``rainbow"
subgraphs (containing all the colors used for the edge-coloring).
This is stated by the next proposition, which is proved in Section
\ref{sec-3}:
\begin{proposition}\label{prop-ramsey}
For every (large) $n$ and $t\leq \sqrt{\frac{2\log n}{(\log\log n)^3
}}$ there is an explicit $t$-edge-coloring of the complete graph on
$n$ vertices, such that every induced subgraph on
$$\mathrm{e}^{(1+o(1))\sqrt{8\log n \log\log n}}$$
vertices contains all $t$ colors.
\end{proposition}
This extends the construction of Frankl and Wilson
\cite{FranklWilson} that deals with the case $t=2$ (using a slightly
different construction).

\section{Graphs with high capacities for unions of predefined
subsets}\label{sec-2} The upper bound on the capacities of subsets
not containing any $F \in \mathcal{F}$ relies on the algebraic bound
for the Shannon capacity using representations by polynomials,
proved in \cite{NogaUnion}. See also Haemers \cite{Haemers} for a
related approach.
\begin{definition} Let $\mathbb{K}$ be a field, and let
$\mathcal{H}$ be a linear subspace of polynomials in $r$ variables
over $\mathbb{K}$. A \textbf{representation} of a graph $G=(V,E)$
over $\mathcal{H}$ is an assignment of a polynomial $f_v \in
\mathcal{H}$ and a value $c_v \in \mathbb{K}^r$ to every $v \in V$,
such that the following holds: for every $v\in V$, $f_v(c_v)\neq 0$,
and for every $u\neq v\in V$ such that $u v \not\in E$,
$f_u(c_v)=0$.
\end{definition}
\begin{theorem}[\cite{NogaUnion}]\label{thm-2.1} Let $G=(V,E)$ be a graph and let
$\mathcal{H}$ be a space of polynomials in $r$ variables over a
field $\mathbb{K}$. If $G$ has a representation over $\mathcal{H}$
then $c(G)\leq \dim(\mathcal{H})$.
\end{theorem}

We need the following simple lemma:
\begin{lemma}\label{lemma-ai-sets} Let $T=[t]$ for $t\geq
1$, and let $\mathcal{F}$ be a family of subsets of $T$. There exist
sets $A_1,A_2,\ldots,A_t$ such that for every $X \subseteq T$:
$$ X \mbox{ does not contain any }F \in \mathcal{F}
\Longleftrightarrow~ \bigcap_{i\in X}A_i \neq \emptyset~.$$
Furthermore, $|\bigcup_{i=1}^t A_i| \leq \binom{t}{\lfloor t/2
\rfloor}$.
\end{lemma}
\begin{proof}[Proof of lemma] Let $\mathcal{Y}$ denote the family of
all maximal sets $Y$ such that $Y$ does not contain any $F \in
\mathcal{F}$. Assign a unique element $p_Y$ to every $Y \in
\mathcal{Y}$, and define: \begin{equation}\label{eq-ai-definition}
A_i = \{ p_Y ~:~i\in Y~,~Y\in\mathcal{Y}\}~.
\end{equation} Let $X
\subseteq T$, and note that \eqref{eq-ai-definition} implies that
$\bigcap_{i\in X}A_i = \{ p_Y ~:~X \subseteq Y\}$. Thus, if $X$ does
not contain any $F \in \mathcal{F}$, then $X\subseteq Y$ for some
$Y\in\mathcal{Y}$, and hence $p_Y \in \bigcap_{i\in X}A_i$.
Otherwise, $X$ contains some $F \in\mathcal{F}$ and hence is not a
subset of any $Y \in \mathcal{Y}$, implying that $\bigcap_{i\in
X}A_i=\emptyset$.

Finally, observe that $\mathcal{Y}$ is an anti-chain and that
$|\bigcup_{i=1}^t A_i| \leq |\mathcal{Y}|$, hence the bound on
$|\bigcup_{i=1}^t A_i|$ follows from Sperner's Theorem
\cite{Sperner}.
\end{proof}
\begin{proof}[\textbf{\em Proof of Theorem \ref{thm-1}}]
Let $p$ be a large prime, and let $\{p_Y ~:~ Y\in\mathcal{Y}\}$ be
the first $|\mathcal{Y}|$ primes succeeding $p$. Define $s = p^2$
and $r=p^3$, and note that, as $t$ and hence $|\mathcal{Y}|$ are
fixed, by well-known results about the distribution of prime
numbers, $p_Y=(1+o(1))p<s$ for all $Y$, where the $o(1)$-term tends
to 0 as $p\to\infty$.

The graph $G_i=(V_i,E_i)$ is defined as follows: its vertex set
$V_i$ consists of all $\binom{r}{s}$ possible $s$-element subsets of
$[r]$, and for every $A\neq B\in V_i$:
\begin{equation}\label{eq-a-b-adjacent}
(A,B) \in E_i ~\Longleftrightarrow~|A\cap B|\equiv
s\pmod{p_Y}
 \mbox{ for some }p_Y \in A_i~.\end{equation}
Let $X \subseteq T$. If $X$ does not contain any $F\in\mathcal{F}$,
then, by Lemma \ref{lemma-ai-sets}, $\bigcap_{i\in X}A_i \neq
\emptyset$, hence there exists some $q$ such that $q \in A_i$ for
every $i \in X$. Therefore, for every $i\in X$, if $A,B$ are
disconnected in $G_i$, then $|A \cap B| \not\equiv s\pmod{q}$. It
follows that the graph $\sum_{i \in X}G_i$ has a representation over
a subspace of the multi-linear polynomials in $|X|r$ variables over
$\mathbb{Z}_q$ with a degree smaller than $q$. To see this, take the
variables $x_j^{(i)}$, $i=1,\ldots,|X|$, $j=1,\ldots,r$, and assign
the following polynomial to each vertex $A\in V_i$:
$$f_A(\overline{x})=\prod_{u \not\equiv s}(u-\sum_{j\in A}x_j^{(i)})~.$$
The assignment $c_A$ is defined as follows: $x_j^{(i')}=1$ if $i'=i$
and $j \in A$, otherwise $x_j^{(i')}=0$. As every assignment
$c_{A'}$ gives values in $\{0,1\}$ to all $x_j^{(i)}$, it is
possible to reduce every $f_A$ modulo the polynomials
$(x_j^{(i)})^2-x_j^{(i)}$ for all $i$ and $j$, and obtain
multi-linear polynomials, equivalent on all the assignments
$c_{A'}$.

The following holds for all $A \in V_i$:
$$f_A(c_A)=\prod_{u \not\equiv s}(u-s)\not\equiv
0\pmod{q}~,$$ and for every $B \neq A$:
\begin{eqnarray}
B \in V_i~,~(A,B)\not\in E_i &\Longrightarrow& f_A(c_B)= \prod_{u
\not\equiv s}(u-|A\cap B|)\equiv 0 \pmod{q}~,\nonumber\\
B \not\in V_i &\Longrightarrow& f_A(c_B)= \prod_{u \not\equiv s}u
\equiv 0 \pmod{q}~,\nonumber
\end{eqnarray}
where the last equality is by the fact that $s \not\equiv 0
\pmod{q}$, as $s=p^2$ and $p < q$. As the polynomials $f_A$ lie in
the direct sum of $|X|$ copies of the space of multi-linear
polynomials in $r$ variables of degree less than $q$, it follows
from Theorem \ref{thm-2.1} that the Shannon capacity of $\sum_{i\in
X}G_i$ is at most:
$$|X|\sum_{i=0}^{q-1}\binom{r}{i}\leq t\sum_{i=0}^{q-1}\binom{
r}{i} < t \binom{r}{q}~.$$ Recalling that $q=(1+o(1))p$ and writing
$t\binom{r}{q}$ in terms of $n=\binom{r}{s}$ gives the required
upper bound on $c(\sum_{i\in X}G_i)$.

Assume now that $X$ contains some $F\in\mathcal{F}$,
$F=\{i_1,\ldots,i_{|F|}\}$. We claim that the following set is an
independent set in $\left(\sum_{i\in X}G_i\right)^{|F|}$:
$$\{ (A^{(i_1)},A^{(i_2)},\ldots,A^{(i_{|F|})})~:~A\subseteq
[r]~,~|A|=s\}~,$$ where $A^{(i_j)}$ is the vertex corresponding to
$A$ in $V_{i_j}$. Indeed, if $(A,A,\ldots,A)$ and $(B,B,\ldots,B)$
are adjacent, then for every $i\in F$, $|A \cap B|\equiv
s\pmod{p_Y}$ for some $p_Y\in A_i$. However, $\bigcap_{i\in
F}A_i=\emptyset$, hence there exist $p_Y\neq p_Y'$ such that $|A
\cap B|$ is equivalent both to $s\pmod{p_Y}$ and to $s\pmod{p_Y'}$.
By the Chinese Remainder Lemma, it follows that $|A\cap B|=s$ (as
$|A\cap B| < p_Y p_Y'$), thus $A=B$. Therefore, the Shannon capacity
of $\sum_{i\in X}G_i$ is at least $\binom{r}{s}^{1/|F|} =
n^{1/|F|}$.
\end{proof}

\section{Explicit construction for rainbow Ramsey
graphs}\label{sec-3} \begin{proof}[\textbf{\em Proof of Proposition
\ref{prop-ramsey}}] Let $p$ be a large prime, and let
$p_1<\ldots<p_t$ denote the first $t$ primes succeeding $p$. We
define $r,s$ as in the proof of Theorem \ref{thm-1}: $s=p^2$,
$r=p^3$, and consider the complete graph on $n$ vertices, $K_n$,
where $n=\binom{r}{s}$, and each vertex corresponds to an
$s$-element subset of $[r]$. The fact that $t \leq \sqrt{\frac{2
\log n}{(\log\log n)^3}}$ implies that $t \leq
(\frac{1}{2}+o(1))\frac{p}{\log p}$, and hence, by the distribution
of prime numbers, $p_t < 2p$ (with room to spare) for a sufficiently
large value of $p$.

We define an edge-coloring $\gamma$ of $K_n$ by $t$ colors in the
following manner: for every $A,B\in V$, $\gamma(A,B) = i$ if $|A\cap
B|\equiv s\pmod{p_i}$ for some $i\in[t]$, and is arbitrary
otherwise. Note that for every $i\neq j \in \{1,\ldots,t\}$, $s <
p_i p_j$. Hence, if $|A\cap B|\equiv s \pmod{p_i}$ and $|A\cap B|
\equiv s \pmod{p_j}$ for such $i$ and $j$, then by the Chinese
Remainder Lemma, $|A\cap B|=s$, and in particular, $A=B$. Therefore,
the coloring $\gamma$ is well-defined.

It remains to show that every large induced subgraph of $K_n$ has
all $t$ colors according to $\gamma$. Indeed, this follows from the
same consideration used in the proof of Theorem \ref{thm-1}. To see
this, let $G_i$ denote the spanning subgraph of $K_n$ whose edge set
consists of all $(A,B)$ such that $\gamma(A,B)\neq i$. Each such
pair satisfies $|A\cap B|\not\equiv s\pmod{p_i}$, hence $G_i$ has a
representation over the multi-linear polynomials in $r$ variables
over $\mathbb{Z}_{p_i}$ with a degree smaller than $p_i$ (define
$f_A(x_1,\ldots,x_r)$ as is in the proof of Theorem \ref{thm-1}, and
take $c_A$ to be the characteristic vector of $A$). Thus, $c(G_i) <
\binom{r}{p_i}$, and in particular, $\alpha(G_i) < \binom{r}{p_i}$.
This ensures that every induced subgraph on at least $\binom{r}{p_i}
\leq \binom{r}{2p}$ vertices contains an $i$-colored edge, and the
result follows.
\end{proof}
\noindent\textbf{Acknowledgement} We would like to thank Benny
Sudakov for fruitful discussions.

\end{document}